\documentclass[12pt,a4paper]{article}
\usepackage{amssymb,amsmath}
\newtheorem{thm}{Theorem}[section]
\newtheorem{cor}[thm]{Corollary}
\newtheorem{lem}[thm]{Lemma}

\newtheorem{rem}[thm]{Remark}

\usepackage[T1]{fontenc}

\title{Extrinsic curvatures of distributions of arbitrary codimension}
\author{Krzysztof Andrzejewski and Pawe\l\, G. Walczak }
\date{August 10,  2009}
\begin{document}
\def\ZZ{\mathbb Z}
\def\RR{\mathbb R}
\def\calf{\mathcal F}
\def\ric{\operatorname{Ric}}
\def\tr{\operatorname{Tr}}
\newcommand\brac[1]{\langle{#1}|}
\newcommand\ket[1]{|{#1}\rangle}
\newcommand\bracket[2]{\langle{#1},{#2}\rangle}
\newcommand\codim{\operatorname{codim}}
\newcommand\divergence{\operatorname{div}}
\newcommand\vol{\operatorname{vol}}

\maketitle

\begin{abstract}
In this article,  using the generalized Newton transformations, we define higher order mean curvatures
of distributions of arbitrary codimension and we show that they coincide with the ones from Brito and Naveira
(Ann. Global Anal. Geom. {\bf 18}, 371--383 (2000)).
We also introduce  higher order mean curvature vector fields   and we compute their 
divergence for certain distributions and using this we  obtain total extrinsic mean curvatures.
\end{abstract}

\section{Introduction}
Using  some special forms $\Gamma_r$ Brito and Naveira \cite{bn} defined   higher 
order extrinsic curvatures of distributions and they computed the total  $r$th mean curvature 
$S_r^T$ of certain distributions on closed spaces of constant curvature.
They  generalize the ones for  foliations \cite{as,blr,ra,ro,wa}.
\par
On the other hand, many authors (see, among the others,
\cite{alm,bc,cl,gr,ro}) have recently  investigated higher order mean curvatures and higher
 order mean curvature vector fields  of hypersurfaces using the Newton transformations 
of the second fundamental form. Especially, the papers \cite{cl,gr}
are devoted to  submanifolds of  codimension greater than one.
In this paper we show that these methods can be also applied successfully for distributions of arbitrary
codimension. Namely, using the  generalized Newton transformation $T_r$ we define
$r$th mean curvature $S_r$ and $(r+1)$th mean curvature vector field $\boldsymbol S_r$ of
a distribution $D$. We show that they  agree with the ones from \cite{bn} (Theorem \ref{t:1}). 
Since most of the interesting
and useful integral formulae in Riemannian geometry are obtained by
computing the divergence of certain vector fields and applying the
divergence theorem, we  compute the divergence of $(r+1)$th mean curvature vector field
 of a distribution which is orthogonal to a totally geodesic foliation in  a manifold of constant 
sectional curvature (Theorem \ref{t:3}).  Using this quantity we obtain a recurrence formula for the total mean curvatures 
(Corollary \ref{c:1})
and consequently we get another proof of the main theorem from \cite{bn} (Theorem \ref{t:2}). 
\par
The paper is organized as follows. Section \ref{s:1} provides some preliminaries.
The main results of the paper are contained in Section \ref{s:2}.
 Throughout the paper everything (manifolds, distribution, metrics and etc.) is assumed to be
$C^{\infty}$-differentiable and oriented and we usually work with $S_r$ instead of its
normalized counterpart $H_r$.
\section{Preliminaries}
\label{s:1}
Let $M$ be a $m$-dimensional oriented, connected Riemannian  manifold.
On $M$ we consider  a distribution $D$, $n=\dim D$  and a distribution $F$ which 
is the  orthogonal complement of $D$, $l=\dim F=m-n$. We assume that both are orientable and
transversally orientable. Let $\langle \cdot ,\cdot\rangle$  represent a metric on $M$
and $\nabla$ denote the Levi-Civita connection of the metric. Let $\Gamma(D)$ denote the set of 
all  vector fields tangent to $D$.
If $v$ is a vector tangent to $M$, then we write  
\[
v=v^\top+v^\bot,
\]
where $v^\top$ is tangent to $D$  and $v^\bot$  to $F$. 
Define the second fundamental form $B$ of the distribution $D$, by 
\[
B(X,Y)=(\nabla_YX)^\bot,
\]
where $X,Y$ are vector fields tangent to $D$. The second fundamental form
of  $F$ is defined similarly. 
\par Throughout this  paper we will use the following index convention:  $1\leq i,j,\ldots\leq n,$
\quad $n+1\leq\alpha,\beta,\ldots \leq m,$ and $1\leq A,B,\ldots\leq m$. Repeated
indices denote summation over their range.
Let us take a local orthonormal frame $\{e_A\}$ adapted to $D,F$, i.e.,  $\{e_i\}$ are tangent  to $D$  and $\{e_{\alpha}\}$ are tangent to $F$. 
Moreover, the frames $\{e_A\}$,$\{e_i\}$ and $\{e_\alpha\}$ are compatible with the 
orientation of $M,D$  and $F$, respectively.  Let  $\{\theta^i\}$ and $\{\theta^\alpha\}$
be their dual frame and  $\omega^{AB}(e_C)=-\langle\nabla_{e_C}e_A,e_B\rangle$
\par Define the second fundamental form (or the shape operator) $A^\alpha$ of
$D$ with respect to $e_\alpha$, by
\[ 
A^\alpha(X)=-(\nabla_Xe_\alpha)^\top,
\]
for $X$ tangent to $D$. Then, using the notation 
\[
A^\alpha e_i={A^\alpha}_i^je_j\quad {\rm and}\quad B^i_j=B(e_i,e_j),
\]
we have 
\[
B^i_j={A^\alpha}_j^ie_\alpha.
\]
Note that,  matrices ${A^\alpha}_j^i$ and $B_j^i$ are not symmetric with respect to $i,j$
if $D$ is not integrable.
In spite of this, for even $r\in\{1,\ldots,n\}$,  we can define   $r$th mean curvature $S_r$
of the distribution $D$ by
\[
S_r=\frac{1}{r!}\delta^{i_1\ldots i_r}_{j_1\ldots j_r}\langle B_{i_1}^{j_1},B_{i_2}^{j_2}\rangle
\cdots\langle B_{i_{r-1}}^{j_{r-1}},B_{i_r}^{j_r}\rangle,
\]
where the  generalized Kronecker symbol  $\delta^{i_1\ldots i_r}_{j_1\ldots j_r}$ is $+1$ or $-1$
according as the $i$'s are distinct and  the $j$'s are either even or odd permutation of the $i$'s,
and is $0$ in all other cases.
By convention, we put  $S_0=1$ and $S_{n+1}=0$.
\par
Moreover, for even $r\in\{0,\ldots,n-1\}$ we define  $(r+1)$th  mean curvature vector field $\boldsymbol{S}_{r+1}$
of $D$ by
\[
\boldsymbol{S}_{r+1}=\frac{1}{(r+1)!}\delta^{i_1\ldots i_{r+1}}_{j_1\ldots j_{r+1}}
\langle B_{i_1}^{j_1},B_{i_2}^{j_2}\rangle
\cdots\langle B_{i_{r-1}}^{j_{r-1}},B_{i_r}^{j_r}\rangle B_{i_{r+1}}^{j_{r+1}}.
\]
We put $\boldsymbol{S}_{n+1}=0$.
If $D$ is of codimension one, then $\boldsymbol S_{r+1}=S_{r+1}N$  where $N$ is a unit vector field
orthogonal to $D$, see \cite{aw}.
The normalized $r$th mean curvature $H_r$ of a distribution $D$ is defined by
\[
 H_r=S_r\dbinom{n}{r}^{-1}.
\]
Obviously, the functions
$S_r$ ($H_r$ respectively) are  smooth  on  the whole $M$. If the distribution $D$ is 
integrable, then for any point $p\in M$,  $S_r(p)$  coincides with the $r$th mean curvature
at $p$  of the leaf $L$ of foliations which passes through  $p$ \cite{alm,cl} .
\par
Now, we introduce the operators  $T_r:\Gamma(D)\rightarrow \Gamma(D)$ which generalizes  the
Newton transformations of the shape operator for hypersurfaces and foliations 
(see, among the others,  \cite{aw,bc,cl,gr,ro}).
\par
 For even $r\in\{1,\ldots,n\}$, we set
\[
{T_r}^i_j=\frac{1}{r!}\delta^{i_1\ldots i_ri}_{j_1\ldots j_rj}\langle B_{i_1}^{j_1},B_{i_2}^{j_2}\rangle
\cdots\langle B_{i_{r-1}}^{j_{r-1}},B_{i_r}^{j_r}\rangle,
\]
and by  convention $T_0=I$. Note that $T_n=0$. 
We  also set for a fixed index $\alpha$  
\[
{T_{r-1}^{\alpha}}^i_j=\frac{1}{(r-1)!}\delta^{i_1\ldots i_{r-1}i}_{j_1\ldots j_{r-1}j}
\langle B_{i_1}^{j_1},B_{i_2}^{j_2}\rangle
\cdots\langle B_{i_{r-3}}^{j_{r-3}},B_{i_{r-2}}^{j_{r-2}}\rangle{A^\alpha}_{i_{r-1}}^{j_{r-1}}.
\]
\par
In the   following lemma, we provide some relations between the $r$th  mean curvature (vector field) 
and the operator $T_r$. 
\begin{lem}
\label{l:1}
For any even integer $r\in\{1,\ldots,n\}$ we have
\begin{align*}
&S_r=\frac 1r\tr(T_{r-1}^\alpha A^\alpha),\\
&\boldsymbol{S}_{r+1}=\frac{1}{r+1}\tr(T_rA^\alpha)e_{\alpha},\\
&\tr(T_r)=(n-r)S_r,\\
&T_r=S_rI-A^\alpha T_{r-1}^{\alpha},\\
\end{align*}
and when $r$ is odd, for each $\alpha$, we have 
\[
tr(T_r^\alpha)=\frac{n-r}{r}\tr(T_{r-1}A^\alpha),
\]
where $\tr=\tr_D=(\cdot)_i^i$. 
\end{lem}
{\it Proof}. 
The proof of lemma is quite similar to the one for submanifolds \cite{cl,gr}, we must only be more 
careful because $B_i^j $ need not be  a symmetric matrix.\hfill$\square$
\par 
On the other hand, Brito and Naveira \cite{bn} have introduced $n$-forms
$\Gamma_r$ for  even $r=2s$ as follows:
\begin{align}
\label{e:1}
\Gamma_r=&\sum_{\sigma\in \Sigma_n}\varepsilon (\sigma)(\omega^{\sigma(1)\beta_1}
\wedge\omega^{\sigma(2)\beta_1})\wedge\cdots\wedge (\omega^{\sigma(2s-1)
\beta_s}\wedge\omega^{\sigma(2s)\beta_s})\wedge\nonumber\\
&\wedge\theta^{\sigma(2s+1)}\wedge \cdots\wedge \theta^{\sigma(n)},
\end{align}
where $\Sigma_n$ is the group of permutations of the set $\{1,\ldots,n\}$, $\varepsilon(\sigma)$
stands for the sign of the permutation $\sigma$.
Furthermore, they define the total $r$th extrinsic mean curvature $S_r^{T}$ of a distribution $D$ on a compact
manifold $M$ as
\[
S_r^T=\frac{1}{r!(n-r)!}\int_M\Gamma_{r}\wedge \nu,
\]
where $\nu=\theta^{n+1}\wedge\cdots\wedge\theta^m$.
This suggests that  we should have
\[
\frac{1}{r!(n-r)!}\Gamma_r\wedge\nu=S_r\Omega,
\]
where $\Omega$ is volume element of $(M,\langle,\rangle)$.
We will show this equality in the next section.
\section{Main results}
\label{s:2}
Using definitions and notations as in Preliminaries, we obtain the following theorem
which states that,  $S_r^T$ defined by Brito and Naveira \cite{bn} is indeed 
the total mean curvature of the distribution in our sense.
\begin{thm}
\label{t:1}
If $r=2s$, $S_r$ is $r$th mean curvature  of the distribution $D$ and $\Gamma_r$ is defined by
(\ref{e:1}), then  we have
\[
\frac{1}{r!(n-r)!}\Gamma_r\wedge\nu=S_r\Omega.
\]
\end{thm}
{\it Proof}.
Using the following expression for the generalized Kronecker symbol  
\[
\delta_{j_1\ldots j_r}^{i_1\ldots i_r}=\left|
\begin{array}{ccc}
\delta^{i_1}_{j_1}&\cdots&\delta^{i_1}_{j_r}\\
\vdots&\ddots&\vdots\\
\delta^{i_r}_{j_1}&\cdots&\delta^{i_r}_{j_r}
\end{array}
\right|=
\sum_{\tau\in\Sigma_r}\varepsilon(\tau)\delta^{i_1}_{j_{\tau(1)}}\cdots\delta^{i_r}_{j_{\tau(r)}},
\]
we have
\begin{align}
\label{e:2}
S_r=&\frac{1}{r!}\delta_{j_1\ldots j_r}^{i_1\ldots i_r}{A^{\alpha_1}}^{j_1}_{i_1}{A^{\alpha_1}}^{j_2}_{i_2}
\cdots{A^{\alpha_s}}^{j_{2s-1}}_{i_{2s-1}}{A^{\alpha_s}}^{j_{2s}}_{i_{2s}}
\nonumber\\
=&\frac{1}{r!}\sum_{\substack{j_1\ldots j_r\\ \rm{distinct}}}\delta_{j_1\ldots j_r}^{i_1\ldots i_r}{A^{\alpha_1}}^{j_1}_{i_1}{A^{\alpha_1}}^{j_2}_{i_2}\cdots
{A^{\alpha_s}}^{j_{2s-1}}_{i_{2s-1}}{A^{\alpha_s}}^{j_{2s}}_{i_{2s}}
\nonumber\\
=&\frac{1}{r!}\sum_{\substack{j_1\ldots j_r\\ \rm{distinct}}}\sum_{\tau\in\Sigma_r}
\varepsilon(\tau)\delta^{i_1}_{j_{\tau(1)}}\cdots\delta^{i_r}_{j_{\tau(r)}}
{A^{\alpha_1}}^{j_1}_{i_1}{A^{\alpha_1}}^{j_2}_{i_2}
\cdots{A^{\alpha_s}}^{j_{2s-1}}_{i_{2s-1}}{A^{\alpha_s}}^{j_{2s}}_{i_{2s}}
\nonumber\\
=&\frac{1}{r!}\sum_{\substack{j_1\ldots j_{r} \\\rm{distinct}}}\sum_{\tau\in\Sigma_r}
\varepsilon(\tau){A^{\alpha_1}}^{j_1}_{j_{\tau(1)}}{A^{\alpha_1}}^{j_2}_{j_{\tau(2)}}\cdots
{A^{\alpha_s}}^{j_{2s-1}}_{j_{\tau(2s-1)}}{A^{\alpha_s}}^{j_{2s}}_{j_{\tau(2s)}}.
\end{align}
On the other hand, by the definition of  $\omega^{i\alpha}$, we deduce
\[
\omega^{i\alpha}(e_j)=\langle e_i,\nabla_{e_j}e_\alpha\rangle={-A^\alpha}_j^i,
\]
thus 
\begin{equation}
\label{e:3}
\omega^{i\alpha}=-{A^\alpha}_j^i\theta^j+X^{i\alpha}_\beta\theta^\beta.
\end{equation}
From  (\ref{e:1}) and (\ref{e:3}),  we have
\begin{align*}
\Gamma_{r}\wedge \nu =&\sum_{\sigma\in\Sigma_n}
\varepsilon(\sigma)({A^{\alpha_1}}^{\sigma(1)}_{j_1}{A^{\alpha_1}}^{\sigma(2)}_{j_2}\cdots
{A^{\alpha_s}}^{\sigma(2s-1)}_{j_{2s-1}}{A^{\alpha_s}}^{\sigma(2s)}_{j_{2s}}\theta^{j_1}
\wedge \cdots\wedge\theta^{j_{2s}})\wedge
\nonumber\\
&\wedge\theta^{\sigma(2s+1)}\wedge\cdots\wedge
\theta^{\sigma(2s)}\wedge\nu\nonumber\\
=&
\sum_{\sigma\in\Sigma_n}\varepsilon(\sigma)\sum_{\tau\in\Sigma\{\sigma(1)\ldots\sigma(2s)\}}
\left(\varepsilon(\tau) {A^{\alpha_1}}^{\sigma(1)}_{\tau(\sigma(1))}
{A^{\alpha_1}}^{\sigma(2)}_{\tau(\sigma(2))} \cdots\right. \nonumber\\ 
&\left. {A^{\alpha_s}}^{\sigma(2s-1)}_{\tau(\sigma({2s-1))}}{A^{\alpha_s}}^{\sigma(2s)}_{\tau(\sigma(2s))}\right)
\theta^{\sigma(1)}\wedge \cdots\wedge\theta^{\sigma(n)}\wedge\nu\nonumber\\
=&\sum_{\sigma\in\Sigma_n}\Big(
\sum_{\tau\in\Sigma\{\sigma(1)\ldots\sigma(2s)\}}\varepsilon
(\tau) {A^{\alpha_1}}^{\sigma(1)}_{\tau(\sigma(1))}{A^{\alpha_1}}^{\sigma(2)}_{\tau(\sigma(2))} \cdots
\nonumber\\ 
&{A^{\alpha_s}}^{\sigma(2s-1)}_{\tau(\sigma({2s-1))}}
{A^{\alpha_s}}^{\sigma(2s)}_{\tau(\sigma(2s))}\Big)\Omega\nonumber\\
=& (n-2s)!\sum_{\sigma:\{1...2s\}\rightarrow \{1...n\}}\Big(
\sum_{\tau\in\Sigma\{\sigma(1)\ldots\sigma(2s)\}}\varepsilon
(\tau) {A^{\alpha_1}}^{\sigma(1)}_{\tau(\sigma(1))}{A^{\alpha_1}}^{\sigma(2)}_{\tau(\sigma(2))} \cdots
\nonumber\\ 
&{A^{\alpha_s}}^{\sigma(2s-1)}_{\tau(\sigma({2s-1))}}{A^{\alpha_s}}^{\sigma(2s)}_{\tau(\sigma(2s))}\Big)
\Omega\nonumber\\
=& (n-2s)!\sum_{\substack{j_1\ldots j_{2s}\\\rm{distinct}}}\Big(
\sum_{\tau\in\Sigma_{2s}}\varepsilon
(\tau) {A^{\alpha_1}}^{j_1}_{j_{\tau(1)}}{A^{\alpha_1}}^{j_2}_{j_{\tau(2)}} \cdots
\nonumber\\ 
&{A^{\alpha_s}}^{j_{2s-1}}_{j_{\tau(2s-1)}}{A^{\alpha_s}}^{j_{2s}}_{j_{\tau(2s)}}\Big)
\Omega.\nonumber\\
\end{align*}
Comparing the above with (\ref{e:2}) we complete the proof of our theorem.\hfill$\square$
\par
Brito and Naveira  have also shown that in some  special cases one can compute explicitly the total mean curvature
$S_r^T$ of the distribution $D$ and  it does not depend on $D$. Indeed, we have the following 
theorem  \cite{bn}.
\begin{thm}
\label{t:2}
If $M$ is a closed manifold of constant sectional curvature $c\geq 0$   and $F=D^{\bot}$ is
a totally geodesic distribution, then 
\[
S_{2s}^T=\left\{
\begin{array}{l}
\dbinom{n/2}{s}\dbinom{l+2s-1}{2s}\dbinom{(l+2s-1)/2}{s}^{-1}c^s\vol(M)\\
\textrm{ if $n$ is even and $l$ is odd,}\\
2^{2s}(s!)^2((2s)!)^{-1}\dbinom{l/2+s-1}{s}\dbinom{n/2}{s}c^s\vol(M)\\
\textrm{ if $n$ and $l$ are even,}\\
0, \textrm{ otherwise.}
\end{array}
\right.
\]
\end{thm}
\begin{rem}
{\rm
Since the distribution $F$ determines a totally geodesic foliation $\calf$ on $M$, the constant
curvature $c$ must be nonnegative; see \cite{ze}.}
\end{rem}
The next part of this section will be devoted to the calculaction of the divergence of the mean curvature vector field.
Next, we will use this to find  a recurrence formula for the total mean curvatures and consequently 
we will  get an alternative proof of Theorem \ref{t:2}. In order to do this we need the following lemma.
\begin{lem}
\label{l:2}
Let $p \in M$ and $ \{e_1,\ldots,e_m\}$  be a local orthonormal frame field adapted to $D$ and $F$, such that
$(\nabla_Xe_i)^{\top}(p)=0$ and $(\nabla_Xe_\alpha)^{\bot}(p)=0$ for any vector field $X$ on
$M$. Then at the point $p$
\begin{align*}
e_\alpha({A^{\beta}}^i_j)=&(A^\beta A^\alpha)_j^i-\langle R(e_j,e_\alpha)e_i,e_\beta\rangle
\\
+&\langle(\nabla _{e_\alpha} e_\gamma)^{\top},e_j\rangle\langle e_i,(\nabla_{e_\gamma}e_\beta)^\top
\rangle-\langle \nabla_{e_j}(\nabla _{e_\alpha}e_\beta)^{\top},e_i\rangle.
\end{align*}
\end{lem}
{\it Proof}.
Our proof starts with the observation that at $p$ we have the following equality 
\[
0=\langle \nabla_{e_j}\nabla_{e_\alpha}e_\beta,e_i\rangle+\langle e_\beta,
\nabla_{e_j}\nabla_{e_\alpha}e_i\rangle.
\]
Thus, we have also  at $p$
\begin{align*} 
&-e_\alpha({A^{\beta}}^i_j)+(A^\beta A^\alpha)_j^i-\langle R(e_j,e_\alpha)e_i,e_\beta\rangle
\\
 &=
(A^\beta A^\alpha)_j^i-\langle\nabla_{e_j}\nabla_{e_\alpha}e_i,e_\beta\rangle+
\langle \nabla_{[e_j,e_\alpha]}e_i,e_\beta\rangle
\\
&={A^\beta}_k^i{A^\alpha}_j^k -\langle\nabla_{e_j}\nabla_{e_\alpha}e_i,e_\beta\rangle
+\langle\nabla_{e_j}e_\alpha,e_k\rangle\langle\nabla_{e_k}e_i,e_\beta\rangle
-\langle\nabla_{e_\alpha}e_j,e_\gamma\rangle\langle\nabla_{e_\gamma}e_i,e_\beta\rangle
\\
&={A^\beta}_k^i{A^\alpha}_j^k +\langle\nabla_{e_j}\nabla_{e_\alpha}e_\beta,e_i\rangle
+\langle\nabla_{e_j}e_\alpha,e_k\rangle\langle\nabla_{e_k}e_i,e_\beta\rangle
-\langle\nabla_{e_\alpha}e_j,e_\gamma\rangle\langle\nabla_{e_\gamma}e_i,e_\beta\rangle
\\
&=\langle\nabla_{e_j}\nabla_{e_\alpha}e_\beta,e_i\rangle
-\langle\nabla_{e_\alpha}e_j,e_\gamma\rangle\langle\nabla_{e_\gamma}e_i,e_\beta\rangle
\\
&=\langle\nabla_{e_j}\nabla_{e_\alpha}e_\beta,e_i\rangle
-\langle(\nabla_{e_\alpha}e_\gamma)^\top,e_j\rangle\langle(\nabla_{e_\gamma}e_\beta
)^\top,e_i\rangle
\\
&=\langle\nabla_{e_j}(\nabla_{e_\alpha}e_\beta)^\top,e_i\rangle
-\langle(\nabla_{e_\alpha}e_\gamma)^\top,e_j\rangle\langle(\nabla_{e_\gamma}e_\beta
)^\top,e_i\rangle.
\end{align*}
This ends the proof.\hfill $\square$
\begin{rem}
 {\rm Note that, using parallel transport in $D$ and $F$ respectively, we can always construct  the frame field from Lemma \ref{l:2}. }
\end{rem}
Now, for even $r$,  we introduce auxiliary notations  as follows 
\[
{T_r}^{i_{r+1}i_{r+2}}_{j_{r+1}j_{r+2}}=\frac{1}{r!}\delta^{i_1\ldots i_{r+2}}_{j_1\ldots j_{r+2}}
\langle B_{i_1}^{j_1},B_{i_2}^{j_2}\rangle\cdots\langle B_{i_{r-1}}^{j_{r-1}},B_{i_r}^{j_r}\rangle,
\]
\[
{T_r}^{i_{r+1}i_{r+2}i}_{j_{r+1}j_{r+2}j}=\frac{1}{r!}\delta^{i_1\ldots i_{r+2}i}_{j_1\ldots j_{r+2}j}
\langle B_{i_1}^{j_1},B_{i_2}^{j_2}\rangle\cdots\langle B_{i_{r-1}}^{j_{r-1}},B_{i_r}^{j_r}\rangle.
\]
\begin{lem}
\label{l:3}
\[
{T_r}^{i_{r+1}i_{r+2}}_{j_{r+1}j_{r+2}}=\delta_{j_{r+2}}^{i_{r+2}}{T_r}_{j_{r+1}}^{i_{r+1}}-
\delta_{j_{r+1}}^{i_{r+2}}{T_r}_{j_{r+2}}^{i_{r+1}}-\frac{1}{r-1}{T_{r-2}}^{i_{r-1}i_{r+1}i_r}_{j_{r-1}
j_{r+1}j_{r+2}}{A^\alpha}^{j_{r-1}}_{i_{r-1}}{A^\alpha}^{i_{r+2}}_{i_r}.
\]
\end{lem}
{\it Proof}. 
The proof is analogous  to the one for submanifolds \cite{cl}.\hfill $\square$
\par Now, we are ready to find the divergence of $\boldsymbol {S}_{r+1}$.
\begin{thm}
\label{t:3}
Let $D$ be a distribution on   a Riemannian manifold $M$ with
constant sectional curvature $c$ and
$S_{r}$($\boldsymbol{S}_{r+1}$) its $r$th mean curvature (vector field), for even  $r\in\{0,1,\ldots,n\}$.
Assume that  $F$  is a totally geodesic distribution(equivalently a totally geodesic foliation) 
orthogonal to $D$. Then 
\[
\divergence(\boldsymbol{S}_{r+1})= -(r+2)S_{r+2}+\frac{c(n-r)(l+r)}{r+1}S_r,
\]
where $n=\dim D$, $l=\dim F$.
\end{thm}
{\it Proof}.
Let $\{e_1,\ldots,e_m\}$ be a frame in the neighbourhood of a point $p$  as in Lemma \ref{l:2}.
By  Lemma \ref{l:1}, we have at $p$
\begin{align}
\label{e:4}
\divergence(\boldsymbol{S}_{r+1})&= \frac{1}{r+1}\langle\nabla_{e_i}(\tr(T_rA^\alpha)
e_\alpha),e_i\rangle+\frac{1}{r+1}\langle\nabla_{e_\beta}(\tr(T_rA^\alpha)
e_\alpha),e_\beta\rangle\nonumber\\
&=\frac{1}{r+1}\tr(T_rA^\alpha)\langle\nabla_{e_i}e_\alpha,e_i\rangle+\frac{1}{r+1}
e_\alpha(\tr(T_rA^\alpha))\nonumber\\
&=-\frac{1}{r+1}\tr(T_rA^\alpha)\tr(A^\alpha)+\frac{1}{r+1}e_\alpha(\tr(T_rA^\alpha)).
\end{align}
Using the definition of $T_r$ and  the symmetries  of the generalized Kronecker symbol  we obtain
\begin{align}
\label{e:5}
e_\alpha(\tr(T_rA^\alpha))=&e_\alpha({T_r}^i_j)
{A^\alpha}_i^j+{T_r}^i_je_\alpha({A^\alpha}_i^j)
\nonumber\\
=&\frac{r}{r!}\delta^{i_1\ldots i_ri}_{j_1\ldots j_rj}
\langle B_{i_1}^{j_1},B_{i_2}^{j_2}\rangle\cdots\langle B_{i_{r-3}}^{j_{r-3}},B_{i_{r-2}}^{j_{r-2}}
\rangle {A^\beta}_{i_{r-1}}^{j_{r-1}} e_\alpha({A^\beta}_{i_{r}}^{j_{r}}) {A^\alpha}_{i}^{j}
\nonumber\\
+&{T_r}_j^ie_\alpha({A^\alpha}_i^j)\nonumber\\
=&\frac{1}{r-1}{T_{r-2}}^{i_{r-1}i_ri}_{j_{r-1}j_rj}{A^\beta}_{i_{r-1}}^{j_{r-1}}
e_\alpha({A^\beta}_{i_r}^{j_r}){A^\alpha}_{i}^{j}
+{T_r}_j^ie_\alpha({A^\alpha}_i^j).
\end{align}
Now let us compute the terms on the right hand side of (\ref{e:5}) one by one.
From Lemma \ref{l:2}, under our assumption $(\nabla_{e_\alpha}e_\beta)^\top=0$, we obtain 
\begin{equation}
\label{e:6}
e_\alpha({A^\beta}^{j_r}_{i_r})=(A^\beta A^\alpha)_{i_r}^{j_r}+c\delta_\alpha^\beta
\delta^{j_r}_{i_r}.
\end{equation}
Using  Lemma  \ref{l:1},  Lemma \ref{l:3}  and (\ref{e:6}), we see that  the first term on the right hand side  of (\ref{e:5})
is of the form
\begin{align*}
&\frac{1}{r-1}{T_{r-2}}^{i_{r-1}i_ri}_{j_{r-1}j_rj}{A^\beta}_{i_{r-1}}^{j_{r-1}}
e_\alpha({A^\beta}_{i_r}^{j_r}){A^\alpha}_{i}^{j}
\\
&=\frac{1}{r-1}{T_{r-2}}^{i_{r-1}i_ri}_{j_{r-1}j_rj}{A^\beta}_{i_{r-1}}^{j_{r-1}}
{A^\alpha}_{i}^{j}{A^\beta}^{j_r}_k {A^\alpha}_{i_r}^k+
\frac{c}{r-1}{T_{r-2}}^{i_{r-1}i_ri}_{j_{r-1}j_rj}{A^\beta}_{i_{r-1}}^{j_{r-1}}
{A^\alpha}_{i}^{j}\delta_\alpha^\beta\delta^{j_r}_{i_r}
\\
&=\frac{1}{r-1}{T_{r-2}}^{i_{r-1}i_ri}_{j_{r-1}j_rj}{A^\alpha}_{i_r}^{j_r}{A^\alpha}_{i}^k
{A^\beta}_{i_{r-1}}^{j_{r-1}}  {A^\beta}^{j}_k+
\frac{c}{r-1}{T_{r-2}}^{i_{r-1}ii_r}_{j_{r-1}ji_r}{A^\alpha}_{i_{r-1}}^{j_{r-1}}
{A^\alpha}_{i}^{j}
 \\
&=\frac{1}{r-1}{T_{r-2}}^{i_{r}i_{r-1}i}_{j_{r}j_{r-1}j}{A^\alpha}_{i_r}^{j_r}{A^\alpha}_{i}^k
{A^\beta}_{i_{r-1}}^{j_{r-1}}  {A^\beta}^{j}_k+cr\tr(T_r)
 \\
&=\left(-{T_r}^{i_{r-1}k}_{j_{r-1}j}+\delta_j^k{T_r}_{j_{r-1}}^{i_{r-1}}-\delta_{j_{r-1}}^k{T_r}_j^{i_{r-1}}
\right){A^\beta}^{j_{r-1}}_{i_{r-1}}{A^\beta}^{j}_{k}+cr(n-r)S_r
\\
&=-{T_r}^{i_{r+1}k}_{j_{r+1}j}{A^\beta}^{j_{r+1}}_{i_{r+1}}{A^\beta}^{j}_{k}+\tr(T_rA^\beta)
\tr(A^\beta)-\tr(T_rA^\beta A^\beta)+cr(n-r)S_r
\\
&=-(r+1)(r+2)S_{r+2}+\tr(T_rA^\beta) \tr(A^\beta)-\tr(T_rA^\beta A^\beta)+cr(n-r)S_r.
\end{align*}
By the use of (\ref{e:6}) and Lemma \ref{l:1}, we see that  the second term on the right hand side  of (\ref{e:5})
is of the form
\begin{align}
{T_r}_j^ie_\alpha({A^\alpha}_i^j)=&{T_r}_j^i(A^\alpha A^\alpha)^j_i+cl{T_r}_i^i
=\tr(T_rA^\alpha A^\alpha)+cl\tr(T_r)
\nonumber \\
=&\tr(T_rA^\alpha A^\alpha)+cl(n-r)S_r.
\end{align}
Hence (\ref{e:5}) is of the form
\begin{equation} 
\label{e:7}
e_\alpha(\tr(T_rA^\alpha))=-(r+1)(r+2)S_{r+2}+\tr(T_rA^\alpha) \tr(A^\alpha)+c(r+l)(n-r)S_r.
\end{equation}
Inserting (\ref{e:7}) into (\ref{e:4}) we  complete the proof of theorem.\hfill $\square$
\begin{cor}
\label{c:1}
Let $D$ be a distribution on   a closed Riemannian manifold $M$ with
constant sectional curvature $c\geq0$ and ($S_r^T$) $S_{r}$  its (total)  $r$th mean curvature. 
Let us assume that  $F$ is a totally geodesic distribution orthogonal to $D$. Then 
\[
\int_M S_{r+2}=\int_M\frac{c(n-r)(l+r)}{(r+1)(r+2)}S_r,
\]
equivalently
\begin{equation}
\label{e:9}
 S_{r+2}^T=\frac{c(n-r)(l+r)}{(r+1)(r+2)}S_r^T.
\end{equation}
\end{cor}
\medskip
\par
 Finally,  note that, we can use Corollary \ref{c:1} to prove Theorem \ref{t:2}.
\medskip
\\
{\it Proof of Theorem \ref{t:2}}.  For even $n$  using (\ref{e:9}) and induction one
gets $S_r^T$ as in Theorem \ref{t:2}. When $n$ is odd, then $c$ must be zero, because there is
no totally geodesic foliation on a closed Riemannian manifold of constant positive curvature.
Indeed, without loss of generality, we may assume that $M=S^m$. For the existence of foliations the
sphere should have odd dimension. Since $n$ is odd, the foliation should be even dimensional
and should not contain any compact spherical leaf. Otherwise, we might pull back the Euler
class of the foliation to this spherical even dimensional leaf, proving that it has Euler number zero.
On the other hand, totally geodesic foliations on round spheres should have 
spheres as leaves - contradiction. Consequently $c=0$ and using again  (\ref{e:9}), we  complete
the proof of Theorem \ref{t:2}.\hfill$\square$
\medskip
\par
For $r=0$ there are known applications of  Theorem \ref{t:3} in 
different areas of differential geometry, analysis and mathematical physics;  see - for example \cite{bs2,bw,sv}.
The reader is warmly invited to find them for  other $r$.
\medskip
\\
{\bf Acknowledgement}
\\
We are grateful to  Fabiano Brito  for helpful e-mail discussion.

Krzysztof Andrzejewski (corresponding author) \\
Institute of Mathematics, Polish Academy of Sciences\\
ul. \'Sniadeckich 8, 00-956 Warszawa, Poland\\
and\\
Department of Theoretical Physics II,
University  of \L \' od\'z\\
ul. Pomorska 149/153, 90 - 236  \L \' od\'z, Poland.\\
e-mail: k-andrzejewski@uni.lodz.pl
\vspace{0.5cm}
\\Pawe\l\, G. Walczak \\
Institute of Mathematics, Polish Academy of Sciences\\
ul. \'Sniadeckich 8, 00-956 Warszawa, Poland\\
and\\
Faculty of Mathematics and Informatics, University of \L \' od\'z\\
ul. Banacha 22, 90-238 \L \' od\'z, Poland\\
 e-mail: pawelwal@math.uni.lodz.pl
\end{document}